\documentclass[11pt]{article}
\usepackage{amsmath, amssymb, amsthm}
\usepackage{graphicx} 

\theoremstyle{plain}

\newtheorem{thm}{Theorem}
\newtheorem{prop}[thm]{Proposition}
\newtheorem{lemma}[thm]{Lemma}
\newtheorem{corollary}[thm]{Corollary}

\theoremstyle{definition}

\begin{document}

\def\Cal#1{{\cal#1}}
\def\<{\langle}
\def\>{\rangle}
\def\what{\widehat}
\def\wtil{\widetilde}
\def\Z{{\mathbb Z}}\def\N{{\mathbb N}} \def\C{{\mathbb C}}
\def\Q{{\mathbb Q}}\def\R{{\mathbb R}}

\def\pro{{\rm prod}}
\def\Tw{{\rm Tw}}
\def\St{{\rm St}}
\def\pos{{\rm pos}}
\def\card{{\rm card}}
\def\An{{\rm An}}
\def\source{{\rm src}}
\def\target{{\rm trg}}
\def\basepts{{\cal P}}
\def\simII{\sim_{\rm II}}

\def\Proof{\paragraph{Proof.}}
\def\Remark{\paragraph{Remark.}}
\def\endproof{\hfill$\bullet$\break\medskip}
\def\noproof{\hfill$\bullet$\break}

\def\notation{\paragraph{Notation.}}
\def\ackn{\paragraph{Acknowledgement.}}

\let\demph\textbf
\let\bk\backslash 
\let\iff\Leftrightarrow
\let\liff\Longleftrightarrow
\let\imply\Rightarrow

%
\def\al{\alpha}                 \def\be{\beta}
\def\ga{\gamma}
\def\sig{\sigma}
\def\ep{\epsilon}               \def\varep{\varepsilon}
%
%

\title{{\bf The solution to a conjecture of Tits on the subgroup generated by the squares of the generators of an Artin group}}

\author{{\bf John Crisp\footnote{The first author gratefully acknowledges the support of a Postdoctoral Grant from the Conseil R\'egional de Bourgogne, France, during the preparation of this work.\hfill\break \break
}, Luis Paris}\\ \\
{\small Laboratoire de Topologie,
Universit\'e de Bourgogne,
UMR 5584 du CNRS,}\\
{\small BP 47 870,
21078 Dijon Cedex, France}\\
{\small(e-mail: \texttt{crisp@topolog.u-bourgogne.fr / lparis@u-bourgogne.fr})}\\}

\date{
{\small Oblatum 21-III-2000 \& 1-XII-2000}\\
{\small Published in \emph{Inventiones Mathematicae} (2001).}\\
{\small Digital Object Identifier (DOI) 10.1007/s002220100138. \hskip5mm\copyright \emph{Springer-Verlag 2001}}\\
}

\maketitle

\begin{abstract}
Let $A$ be an Artin group with standard generating set $\{\sigma_s : s\in S\}$. Tits 
conjectured that the only relations in $A$ amongst the squares of the generators are consequences of the obvious ones, namely that $\sig_s^2$ and $\sig_t^2$ commute whenever $\sig_s$ and $\sig_t$ commute, for $s,t\in S$. 
In this paper we prove Tits' conjecture for all Artin groups.
In fact, given a number $m_s\geq 2$ for each $s\in S$, we show that the elements $\{ T_s=\sig_s^{m_s} : s\in S\}$ generate a subgroup that has a finite presentation in which the only defining relations are that $T_s$ and $T_t$ commute if $\sig_s$ and $\sig_t$ commute.
\end{abstract}

\section{Introduction} 

Let $S$ be a finite set. A \demph{Coxeter matrix} over $S$ is a matrix $M=(m_{s,t})_{s,t\in S}$ indexed by the elements of $S$ and such that:
\smallskip

\noindent (a) {$m_{s,s}=1$ for all $s\in S$, and}

\noindent (b) {$m_{s,t}=m_{t,s}\in \{ m\in\N : m\geq 2\}\cup\{\infty\}$ for all $s\neq t\in S$.}

\smallskip

\noindent A Coxeter matrix is represented by its \demph{Coxeter graph} $\Gamma=\Gamma(M)$, which is defined as follows:

\smallskip

\noindent (a) {$S$ is the set of vertices of $\Gamma$;}

\noindent (b) {two vertices $s$ and $t$ are joined by an edge if $s\neq t$ and $m_{s,t}\neq 2$;}

\noindent (c) {the edge joining $s$ and $t$ is labelled $m_{s,t}$ if $m_{s,t}\geq 4$ or $m_{s,t}=\infty$.}

 \smallskip

Take an abstract set $\Sigma = \{ \sig_s : s\in S\}$ in bijection with $S$. For $s,t\in S$ and $m\in\N$ define the word
\[
\pro(\sig_s,\sig_t; m) =
\begin{cases}
(\sig_s\sig_t)^{\frac{m}{2}}\qquad 	&\text{if $m$ is even},\\
(\sig_s\sig_t)^{\frac{m-1}{2}}\sig_s\qquad&\text{if $m$ is odd}.
\end{cases}
\]

\noindent The \demph{Artin system} associated to $\Gamma=\Gamma(M)$ is the pair $(A(\Gamma),\Sigma)$ where $A(\Gamma)$ is the group defined by the presentation
\[
\langle \Sigma\mid \pro(\sig_s,\sig_t;m_{s,t})=\pro(\sig_t,\sig_s;m_{s,t})\,,\text{ for all }s,t\in S \text{ with } m_{s,t}\neq\infty \rangle\,. 
\]

This paper is concerned with subgroups of the Artin group $A=A(\Gamma)$ which are generated by powers of the generators. In \cite{AS}, Appel and Schupp showed that, in the case of so-called ``extra-large type'' Artin groups (i.e: where $m_{s,t}\notin\{ 2,3\}$ for $s,t\in S$), the set of squares of generators $\Cal Q=\{Q_s=\sig_s^2 : s\in S\}$ freely generates a free subgroup of $A$. They then stated the conjecture, communicated to them by Tits, that for an arbitrary Artin group the obvious relations, $Q_sQ_t=Q_tQ_s$ if $m_{s,t}=2$, might be sufficient to give a presentation of the subgroup generated by $\Cal Q$.
This conjecture was apparently motivated by recalling a much earlier result of Tits \cite[Th\'eor\`eme 2.5]{Tits} which says, essentially, that the abelianization of the pure Artin group maps the squares of the generators of  $A$ to a set of independent generators of a free abelian group, showing that the  subgroup in question lies somewhere between the group  $H=\langle\Cal Q\mid Q_sQ_t=Q_tQ_s\text{ if }m_{s,t}=2\rangle$ and the abelianization of $H$.

In this paper we prove the following theorem:

\begin{thm}\label{main}
Let $\Gamma$ be a Coxeter graph over $S$, and suppose that one is given an integer $m_s\geq 2$ for each $s\in S$. Then the set of elements $\Cal T= \{T_s=\sig_s^{m_s} : s\in S\}$ generates a subgroup of $A(\Gamma)$ with presentation
\[
\langle\Cal T\mid T_sT_t=T_tT_s\text{ if }m_{s,t}=2\rangle\,.
\]  
\end{thm}

\begin{corollary}
The Tits Conjecture holds for every Artin system. Namely, given an arbitrary Coxeter graph $\Gamma$, the set of squares of the Artin generators, $\Cal Q=\{Q_s=\sig_s^2 : s\in S\}$, generates a subgroup of $A(\Gamma)$ with presentation
\[
\langle\Cal Q\mid Q_sQ_t=Q_tQ_s\text{ if }m_{s,t}=2\rangle\,.
\]  

\end{corollary}

After the paper of Appel and Schupp, the Tits Conjecture was proven in various special cases: by Pride \cite{Pride} for the ``triangle-free'' groups (i.e: those in which every three generators include a pair $s,t$ for which $m_{s,t}=\infty$); by Droms, Lewin and Servatius \cite{DLS} for braid groups on fewer than $6$ strings; by Collins \cite{Coll} for all braid groups; by Humphries \cite{Hum} in some restricted cases where $m_{s,t}\in \{ 2,3\}$ for distinct $s,t\in S$; and recently by Charney \cite{Ch} for ``locally reducible'' Artin groups (those in which every finite type special subgroup is a direct product of special subgroups of rank at most $2$) which includes the triangle-free and extra-large cases, but none of the irreducible finite types of rank greater than $2$. The arguments given here are self-contained and independent of these previous approaches to the problem. 

We note that there is an obvious homomorphism $f$ from the abstract group $H(\Gamma)$ defined by the presentation displayed in Theorem \ref{main} to the group $A(\Gamma)$. Our task is to show that $f$ is an isomorphism onto its image.
We consider first the Artin groups of \demph{small type}, those for which $m_{s,t}\in\{ 2,3\}$ if $s,t\in S$ are distinct. In Section \ref{S3} we describe, for each irreducible Artin group $A(\Gamma)$ of small type, a representation of $A(\Gamma)$ in the mapping class group of a connected surface $F(\Gamma)$ with boundary.
This naturally induces an action of $A(\Gamma)$, and via the map $f$ an action of $H(\Gamma)$,  on the fundamental groupoid of $F(\Gamma)$ relative to a certain finite set of ``basepoints''.
The result follows by showing that the action of $H(\Gamma)$ on this groupoid is faithful. This is done in Sections \ref{S4} and \ref{S5}. 
The proof proceeds by induction on the length of so-called $M$-reduced expressions in the generators of $H(\Gamma)$, and relies also on the fact that the fundamental groupoid which we consider is isomorphic to a free groupoid on a graph.  Notions concerning groupoids are discussed in Section \ref{S2}. 

We explain in Section \ref{S6} how one may construct a homomorphism from an arbitrary Artin group $A(\Gamma)$ to an Artin group $A(\wtil\Gamma)$ of small type such that the induced map $H(\Gamma)\to H(\wtil\Gamma)$  is injective. This reduces the general case of Theorem \ref{main}
to the small type case already discussed.

\section{The free groupoid on a graph}\label{S2} 

A \demph{groupoid} is a category in which every morphism is invertible and which is connected in the sense that there exists a morphism  between any pair of objects. A group, for example, is a groupoid with exactly one object. The objects of a groupoid are called \demph{vertices}. For every vertex $x$, the morphisms from $x$ to itself form a group called the \demph{vertex group}. An element is said to be \demph{constant} if it is the identity of some vertex group, and \demph{non-constant} otherwise. 

To give another important example of a groupoid, let $X$ be a path-connected topological space and let $\{P_i\}_{i\in I}$ be a discrete set of points in $X$. The \demph{fundamental groupoid} of $X$ based at  $\{P_i\}_{i\in I}$ is the groupoid, denoted $\pi_1(X, \{P_i\}_{i\in I})$, whose vertex set is $\{P_i\}_{i\in I}$, and in which the morphisms from $P_i$ to $P_j$ are the homotopy classes of paths from $P_i$ to $P_j$. Constant morphisms are represented by constant paths.

We now give a combinatorial construction which is motivated by considering  fundamental groupoids of graphs (1-dimensional simplicial complexes).

An \demph{oriented graph} $G$ consists of the following:
\begin{description}
\item{(a)} {a set $V(G)$ of \demph{vertices};}
\item{(b)} {a set $E(G)$ of \demph{edges};}
\item{(c)} {a \demph{source} function $\source :E(G)\to V(G)$ and a \demph{target} function $\target :E(G)\to V(G)$.}
\end{description}

We define an abstract set $E(G)^{-1}=\{ a^{-1}: a\in E(G)\}$ in bijection with $E(G)$ and write $\source(a^{-1})=\target(a)$ and $\target(a^{-1})=\source(a)$ for all $a\in E(G)$. A \demph{path} in $G$ is an expression 
\[
f=a_1^{\varep_1}a_2^{\varep_2}...a_k^{\varep_k}
\]
where $a_i\in E(G)$, $\varep_i\in\{\pm 1\}$, and $\target(a_i^{\varep_i})=\source(a_{i+1}^{\varep_{i+1}})$.
The vertex $\source(a_1^{\varep_1})$ is called the \demph{source} of $g$, written $\source(g)$, and $\target(a_k^{\varep_k})$ the \demph{target} of $g$, written $\target(g)$. The number $k$ is the \demph{length} of $g$. The vertices are considered as paths of length $0$. The \demph{inverse} of a path  
$f=a_1^{\varep_1}a_2^{\varep_2}...a_k^{\varep_k}$ is the path
$f^{-1}=a_k^{-\varep_k}...a_2^{-\varep_2}a_1^{-\varep_1}$. The \demph{product} of two paths $f=a_1^{\varep_1}a_2^{\varep_2}...a_p^{\varep_p}$ and $g=b_1^{\mu_1}b_2^{\mu_2}...b_q^{\mu_q}$ such that $\target(f)=\source(g)$ 
is the path $fg= a_1^{\varep_1}a_2^{\varep_2}...a_p^{\varep_p} b_1^{\mu_1}b_2^{\mu_2}...b_q^{\mu_q}$. A path $f=a_1^{\varep_1}a_2^{\varep_2}...a_k^{\varep_k}$
is \demph{reduced} if there is no $i\in\{ 1,2,..,(k-1)\}$ such that $a_i=a_{i+1}$ and $\varep_i=-\varep_{i+1}$.

A \demph{homotopy} on $G$ is an equivalence relation, $\sim$, on the paths of $G$ satisfying:
\begin{description}
\item{(a)} {$f\sim g \imply \source(f)=\source(g) \text{ and } \target(f)=\target(g)$;}
\item{(b)} {$f^{-1}f\sim\source(f)$  for all paths $f$;}
\item{(c)} {$f\sim g \imply f^{-1}\sim g^{-1}$;}
\item{(d)} {if $f\sim g$, $\target(h_1)=\source(f)=\source(g)$ and $\source(h_2)=\target(f)=\target(g)$, then $h_1fh_2\sim h_1gh_2$. }
\end{description}

A homotopy, $\sim$, on $G$ determines a groupoid $\pi(G,\sim)$ whose objects are the vertices of $G$ and  whose morphisms are the equivalence classes of paths. The \demph{free groupoid} on $G$, denoted $\pi(G)$, is the groupoid determined by the smallest homotopy of $G$. Note that the morphisms of $\pi(G)$ are in bijection with the reduced paths. Since it is the usual convention to compose morphisms in a category on the left, one should associate the domain and codomain of a homotopy class of paths with the target and source of a representative, in that order. (This last point is purely a matter of convention, and of no significance in what follows).  We will use, later, the fact that $\pi(G)$ is isomorphic to the fundamental groupoid $\pi_1(G,V(G))$ via the map defined by sending each edge $e$ to the homotopy class $[\gamma_e]$ where $\gamma_e$ denotes a path in $G$ which travels along the edge $e$ from $\source(e)$ to $\target(e)$.

\section{Representation in the Mapping Class Group of a surface}\label{S3}

Let $F$ be a compact, connected, oriented surface with boundary ($\partial F\neq \emptyset$). We denote by $\Cal H(F)$ the group of homeomorphisms $h:F\to F$ which preserve the orientation and restrict to the identity on the boundary of $F$. The \demph{mapping class group} of $F$ is the group $\Cal M(F)$ of isotopy classes in $\Cal H(F)$.

If $\{P_i\}_{i\in I}$ is a finite collection of points on the boundary of $F$, then $\Cal M(F)$ acts on the fundamental groupoid $\pi_1(F, \{P_i\}_{i\in I})$ by automorphisms.

An \demph{essential circle} in $F$ is an oriented embedding $a:S^1\to F$ such that $a(S^1)\cap\partial F=\emptyset$. We shall use the description of the circle $S^1$ as $\R/2\pi\Z$. Take an (oriented) embedding $A: S^1\times [0,1]\to F$ of the annulus such that $A(\theta,\frac{1}{2})=a(\theta)$ for all $\theta\in S^1$, and define a homeomorphism $\Tw_a\in\Cal H(F)$ which restricts to the identity outside the interior of the image of $A$, and is otherwise given by

\[
(\Tw_a\circ A)(\theta,x) = A(\theta + 2\pi x,x)\quad \text{ for } (\theta,x)\in S^1\times [0,1]\,.
\]

\begin{figure}[ht]
\begin{center}
\includegraphics[width=15cm]{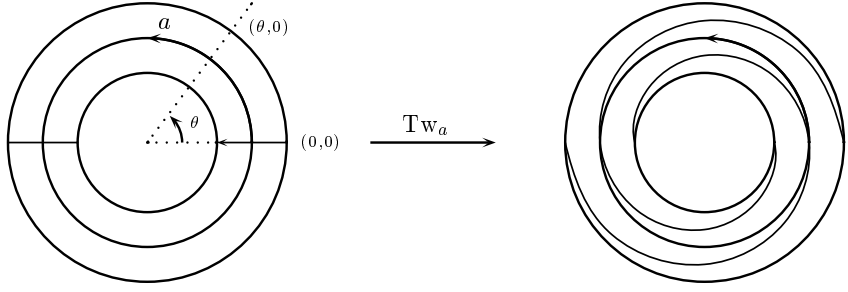}
\end{center}
\caption{A Dehn twist}\label{Fig1}
\end{figure}

\noindent The \demph{Dehn twist} along $a$ is the element of $\Cal M(F)$ represented by $\Tw_a$ (see Figure \ref{Fig1}).
The following result is easily checked and may be found, for instance, in \cite{Bir}. 

\begin{prop}\label{P1}
Let $a_1,a_2:S^1\to F$ be two essential circles, and for $i=1,2$, let $\tau_i$ denote the Dehn twists along $a_i$. Then 
\[\begin{aligned}
\tau_1\tau_2&=\tau_2\tau_1\qquad &\text{if }\ a_1\cap a_2 =\emptyset\\
\tau_1\tau_2\tau_1&=\tau_2\tau_1\tau_2\qquad &\text{if }\ |a_1\cap a_2| =1\,.
\end{aligned}
\]\end{prop}

We suppose now that $\Gamma$ is a connected Coxeter graph of \demph{small type}, which is to say that $m_{s,t}\in\{ 2,3\}$ for distinct $s,t\in S$. We are first going to associate to $\Gamma$ a surface $F(\Gamma)$.

Let $<$ denote a total order on $S$, which may be chosen arbitrarily. For each $s\in S$ we define the set (``star of $s$'')  
\[
\St_s=\{t\in S: m_{s,t}=3\}\cup\{ s\}\,.
\] 

We write $\St_s=\{ t_1,t_2,...,t_k\}$ such that $t_1<t_2<...<t_k$ and suppose that $s=t_j$. The difference $i-j$ will be called the \demph{relative position} of $t_i$ with respect to $s$ and will be denoted $\pos(t_i;s)$. In particular, $\pos(s;s)=0$.

Let $s\in S$. Put $k=\card(\St_s)$.  We denote by $\An_s$ the annulus defined by
\[
\An_s = \R/2k\Z\times[0,1]\ ,
\]
\noindent 
For each $s\in S$, write $P_s$ for the point $(0,0)$ of $\An_s$.

\begin{figure}[ht]
\begin{center}
\includegraphics[width=15cm]{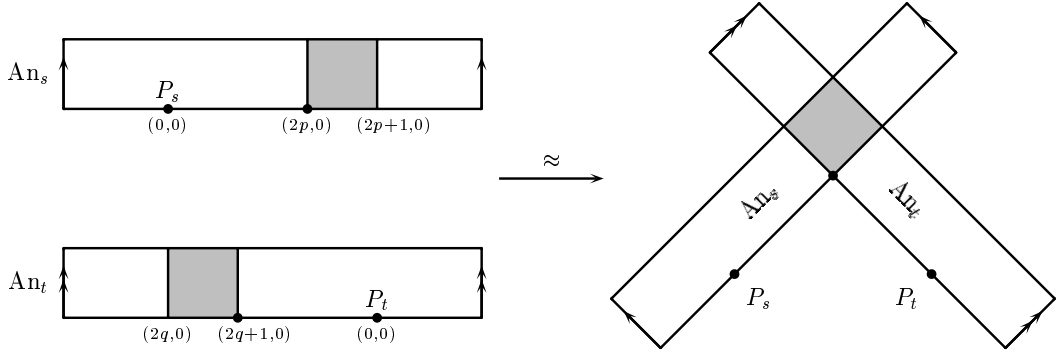}
\end{center}
\caption{Identification of annuli}\label{Fig2}
\end{figure}

The surface $F(\Gamma)$ is defined by
\[
F(\Gamma) = 
\raise.75ex\hbox{$(\coprod\limits_{s\in S}An_s)$}/\lower.75ex\hbox{$\approx$}\ ,
\]
where $\approx$ is the relation defined as follows. Let $s,t\in S$ such that $m_{s,t}=3$ and $s<t$. Put $p=\pos(t;s)>0$ and $q=\pos(s;t)<0$. For each $(x,y)\in [0,1]\times [0,1]$, the relation $\approx$ identifies the point $(2p+x,y)$ of $\An_s$ with the point $(2q+1-y,x)$ of $\An_t$ (see Figure \ref{Fig2}). We identify each annulus $\An_s$ and the point $P_s$ with their images in $F(\Gamma)$ respectively. We also define the set $\basepts = \{ P_s\}_{s\in S}$ in the boundary  of $F(\Gamma)$ to be used as the set of basepoints for a fundamental groupoid.
\medskip

\begin{figure}[ht]
\begin{center}
\includegraphics[width=15cm]{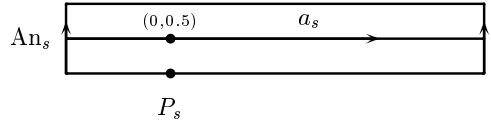}
\end{center}
\caption{The essential circle $a_s$ in $\An_s$}\label{Fig3}
\end{figure}

We now define a group homomorphism $g:A(\Gamma)\to\Cal M(F(\Gamma))$. Let $s\in S$, and put $k=\card(\St_s)$. We denote by $a_s:S^1\to F(\Gamma)$ the essential circle of $F(\Gamma)$ such that $a_s(\theta)$ is the point $(\frac{k\theta}{\pi},\frac{1}{2})$ of $\An_s$ (see Figure \ref{Fig3}). 

We let $\tau_s$ denote the Dehn twist along $a_s$. One has $a_s\cap a_t=\emptyset$ if $m_{s,t}=2$, and $|a_s\cap a_t|=1$ if $m_{s,t}=3$. Therefore, by Proposition \ref{P1}, we have the following.

\begin{prop}\label{P2} Let $\Gamma$ be a connected Coxeter graph of small type. There exists a well-defined group homomorphism $g: A(\Gamma)\to\Cal M(F(\Gamma))$ which sends $\sig_s$ to $\tau_s$ for each $s\in S$. In particular, the group $A(\Gamma)$ acts (via the homomorphism $g$ and the natural action of the mapping class group) by automorphisms of the fundamental groupoid
$\pi_1(F(\Gamma),\basepts)$. 
\end{prop}

We finish this Section with a description of the fundamental groupoid $\pi_1(F(\Gamma),\basepts)$ as the free groupoid of a certain graph.
One associates to $\Gamma$ an oriented graph $G(\Gamma)$ as follows. First let \[\Cal B=\{(s,t)\in S\times S : m_{s,t}=3 \text{ and } s<t\}\,.\] Then put
\[
\begin{aligned}
V(G(\Gamma))&= S\,, \text{ and }\\
E(G(\Gamma))&=\{ e_s: s\in S\}\cup\{ f_{s,t}:(s,t)\in\Cal B\}
\end{aligned}
\]
where $\source(e_s)=\target(e_s)=s$, for all $s\in S$, and $\source(f_{s,t})=s$ and $\target(f_{s,t})=t$, for all $(s,t)\in\Cal B$.

\begin{figure}[ht]\begin{center}
\includegraphics[width=15cm]{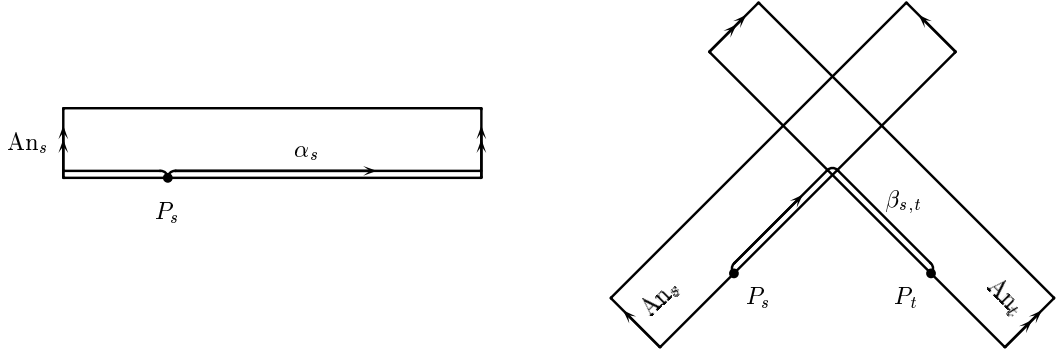}
\end{center}
\caption{The loop $\al_s$ and the path $\be_{s,t}$}\label{Fig4}
\end{figure}

We turn now to the fundamental groupoid.
Fix $s\in S$ and let $k=\card(\St_s)$. Then we define $\al_s:[0,1]\to F(\Gamma)$ to be the loop at $P_s$ such that $\al_s(x)$ is the point $(2kx,0)$ of $\An_s$ (see Figure \ref{Fig4}).

Fix $(s,t)\in\Cal B$ and let $p=\pos(t;s)>0$ and $q=\pos(s;t)<0$. We define $\be_{s,\bullet}:[0,1]\to F(\Gamma)$ to be the path such that $\be_{s,\bullet}(x)$ is the point $(2px,0)$ of $\An_s$ and $\be_{\bullet,t}:[0,1]\to F(\Gamma)$ the path such that $\be_{\bullet,t}(x)$ is the point $((2q+1)(1-x),0)$ of $\An_t$. Then $\be_{s,t}:[0,1]\to F(\Gamma)$ is the path from $P_s$ to $P_t$ defined by $\be_{s,t}=\be_{s,\bullet}\be_{\bullet,t}$ (see Figure \ref{Fig4}).

We now have the following proposition.

\begin{prop}\label{P3}
There is an isomorphism $\pi(G(\Gamma))\to \pi_1(F(\Gamma),\basepts)$ which sends $e_s$ to $[\al_s]$ for all $s\in S$, and $f_{s,t}$ to $[\be_{s,t}]$ for all $(s,t)\in\Cal B$.
\end{prop}

\Proof
Recall that, for an oriented graph $G$, the free groupoid $\pi(G)$ is canonically isomorphic to the fundamental groupoid $\pi_1(G,V(G))$. The Proposition now follows from the observation that the map $G(\Gamma)\to F(\Gamma)$ which sends the edge $e_s$ to the loop $\al_s$, for each $s\in S$, and the edge $f_{s,t}$ to the path $\be_{s,t}$, for each $(s,t)\in\Cal B$, is a homotopy equivalence relative to the vertex set $S$ and its image $\{ P_s\}_{s\in S}$.    
\endproof

\section{Representation by automorphisms of a free groupoid}\label{S4}

Suppose still that $\Gamma$ is a connected Coxeter graph of small type i.e: $m_{s,t}\in\{ 2,3\}$ for distinct $s,t\in S$. We turn now to study in detail 
the action of the Artin group $A(\Gamma)$ on the free groupoid $\pi_1(F(\Gamma),\basepts)$. For convenience, we shall use $\al_s$ and $\be_{s,t}$ to denote the groupoid generators represented by the paths of the same name defined in the previous section. We also use $\tau_s$ to denote both the Dehn twist along the midline $a_s$ of $\An_s$ and the groupoid automorphism which this Dehn twist induces.

Let $B=\< \be_{s,t}: (s,t)\in\Cal B\>$ denote the subgroupoid of $\pi_1(F(\Gamma),\basepts)$ generated by the set  $\{ \be_{s,t} : (s,t)\in\Cal B\}$. This subgroupoid can also be described as follows. Let $G_0(\Gamma)$ denote the subgraph of $G(\Gamma)$ having the same vertex set, but including only the edges $f_{s,t}$ for $(s,t)\in\Cal B$. Then one has an isomorphism $\pi(G_0(\Gamma))\to B$ sending $f_{s,t}$ to $\be_{s,t}$ for all $(s,t)\in\Cal B$.

Observe that any element $w\in\pi_1(F(\Gamma),\basepts)$ has a unique expression of the form \[ w=\mu_0\al_{s_1}^{k_1}\mu_1...\al_{s_\ell}^{k_\ell}\mu_\ell \] where $\mu_i\in B$, $\mu_i$ is non-constant if $i\notin\{ 0,\ell\}$, $s_i\in S$, and $k_i\in\Z\setminus\{ 0\}$. We call this expression the
\demph{reduced form} of $w$.

Let  $s\in S$. Given the reduced form as above, we shall say that $w$ \demph{has a square in $\al_s$} if there exists $i\in\{ 1,..,\ell\}$ such that $s_i=s$ and $|k_i|\geq 2$. Otherwise, we say that $w$ \demph{is without squares in $\al_s$}.  

Now let $\mu\in B$, $t\in S$, and $m\in\Z\setminus\{ 0\}$. we say that $w$ \demph{is of type} $(\mu, \al_t^m)$ if $w$ has reduced form \[
w=\mu_0\al_{t}^{k_1m}\mu_1...\al_{t}^{k_\ell m}\mu_\ell \] where  
$\mu_i\in B$, $\mu_i$ is non-constant if $i\notin\{ 0,\ell\}$, $k_i\in\Z\setminus\{ 0\}$, and $\mu=\mu_0\mu_1...\mu_\ell$.
Note that, for fixed $t\in S$ and $m\in\Z\setminus\{ 0\}$, the set $B(\al_t^m)$ of all elements of type $(\mu, \al_t^m)$ (for some $\mu\in B$) forms a subgroupoid of $\pi_1(F(\Gamma),\basepts)$ which admits a homomorphism onto $B$ sending elements of type $(\mu, \al_t^m)$ to $\mu$. 

Recall that, for $t\in S$, the generator $\sigma_t$ of $A(\Gamma)$ acts on $\pi_1(F(\Gamma),\basepts)$ via the Dehn twist $\tau_t$ along the midline $a_t$ of $\An_t$.
As it will be useful later, in Section \ref{S5}, we now consider the action on the groupoid of an arbitrary power of $\tau_t$ (rather than simply $\tau_t$ itself). We have the following:

\begin{prop}\label{P4}
Let $t\in S$ and $m\in\Z\setminus\{0\}$. 
\begin{description}
\item{(i)}{ Let $\mu\in B$. Then $\tau_t^m(\mu)$ is of type $(\mu,\al_t^m)$.}
\item{(ii)}{ Let $s\in S$. If $s=t$ or if $m_{s,t}=2$, then\[ \tau_t^m(\al_s)=\al_s\,.\] If $m_{s,t}=3$, then
\[ \tau_t^m(\al_s)=
\begin{cases}
u\al_s\quad\text{ if }s<t\\
\al_sv\quad\text{ if }t<s\,,
\end{cases}\] where $u$, respectively $v$, is a non-constant element of $\pi_1(F(\Gamma),\basepts)$ of type $(1,\al_t^m)$. In particular, if $|m|\geq 2$, then $\tau_t^m(\al_s)$ has a square in $\al_t$.}
\end{description}
\end{prop}

\Proof
(i) It suffices to consider the case where $\mu=\be_{r,s}$ for some $(r,s)\in\Cal B$. When $a_t\cap\be_{r,s}=\emptyset$,
we have $\tau_t^m(\be_{r,s})=\be_{r,s}$ as required. On the other hand, we have $a_t\cap\be_{r,s}\neq\emptyset$ if and only if $r<t<s$ and either $t\in\St_r$ or $t\in\St_s$. 
(Recall that $\St_s=\{t\in S: m_{s,t}=3\}\cup\{ s\}$, for $s\in S$.)

\begin{figure}
\begin{center}
\includegraphics[width=15cm]{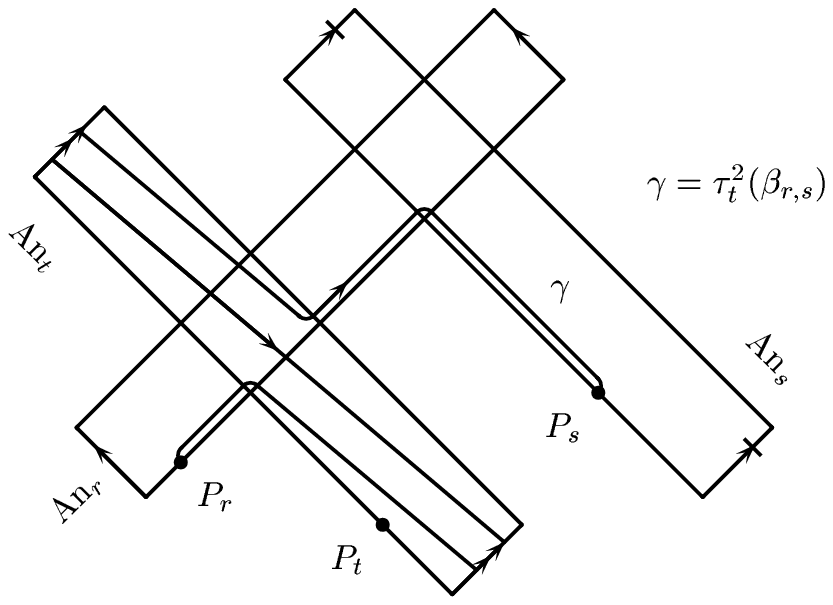}
\end{center}
\caption{The path $\tau^2_t(\be_{r,s})$ when $t\in\St_r$ but $t\not\in\St_s$.}\label{Fig5}
\end{figure}
 
\begin{figure}
\begin{center}
\includegraphics[width=15cm]{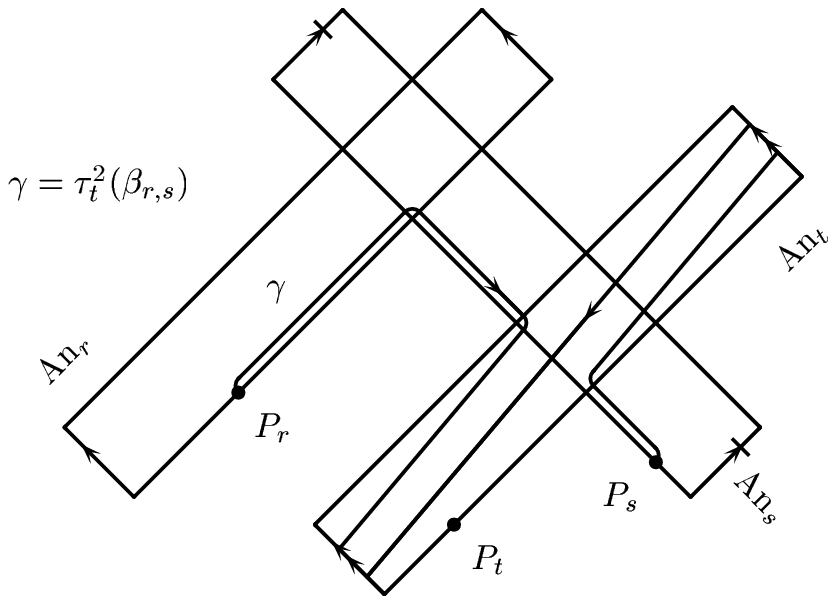}
\end{center}
\caption{The path $\tau^2_t(\be_{r,s})$ when $t\not\in\St_r$ but $t\in\St_s$.}
\label{Fig6}
\end{figure}

\begin{figure}
\begin{center}
\includegraphics[width=15cm]{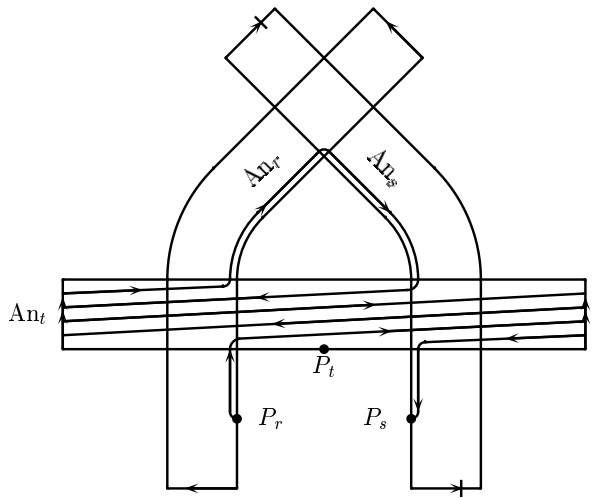}
\end{center}
\caption{The path $\tau^2_t(\be_{r,s})$ when  $t\in\St_s\cap\St_r$.}\label{Fig7}
\end{figure}


Suppose that $t\in\St_r$, but $t\not\in\St_s$. Then

\[
\tau_t^m(\be_{r,s})=\be_{r,t}\al_t^m\be_{r,t}^{-1} \be_{r,s} \qquad\text{ (see Figure \ref{Fig5}).}
\]

Suppose that $t\in\St_s$, but $t\not\in\St_r$.
Then

\[
\tau_t^m(\be_{r,s})=\be_{r,s}\be_{t,s}^{-1}\al_t^{-m}\be_{t,s} \qquad\text{ (see Figure \ref{Fig6}).}
\]

Finally, suppose that $t\in\St_r\cap\St_s$. Then

\[
\tau_t^m(\be_{r,s})=\be_{r,t}\al_t^m\be_{r,t}^{-1} \be_{r,s} \be_{t,s}^{-1}\al_t^{-m}\be_{t,s}\qquad\text{ (see Figure \ref{Fig7}).}
\]

\begin{figure}
\begin{center}
\includegraphics[width=15cm]{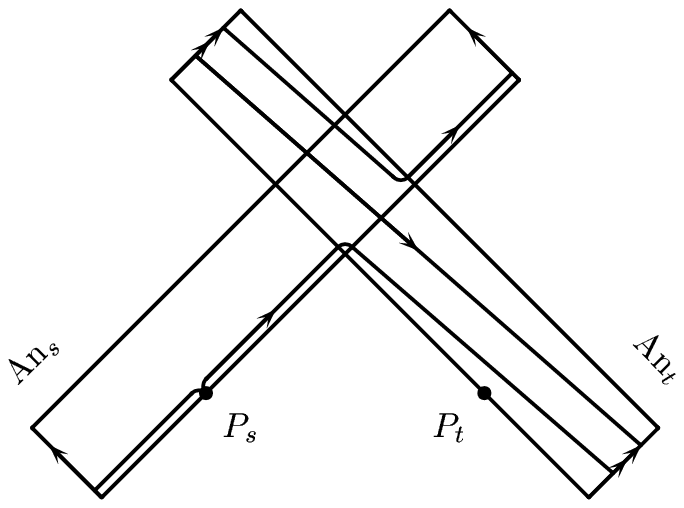}
\end{center}
\caption{The loop $\tau^2_t(\al_s)$ when  $t\in\St_s$ and $s<t$.}\label{Fig8}
\end{figure}

(ii) We have $\al_s\cap a_t\neq\emptyset$ if and only if $m_{s,t}=3$.
In particular, if $s=t$ or if $m_{s,t}=2$, we have $\tau_t^m(\al_s)=\al_s$. 

Suppose that $m_{s,t}=3$. If $s<t$, then
\[
\tau_t^m(\al_s)=\be_{s,t}\al_t^{m}\be_{s,t}^{-1}\al_s \qquad\text{ (see Figure \ref{Fig8}).}
\]

\begin{figure}[ht]
\begin{center}
\includegraphics[width=15cm]{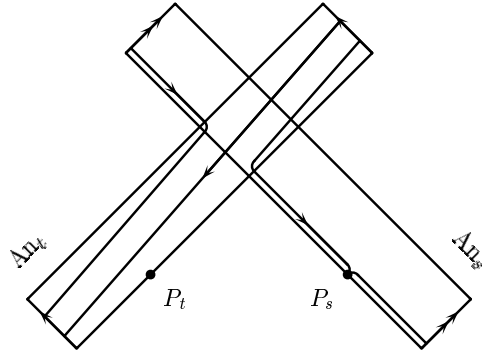}
\end{center}
\caption{The loop $\tau^2_t(\al_s)$ when  $t\in\St_s$ and $s>t$.}\label{Fig9}
\end{figure}

If $s>t$, then

\[
\tau_t^m(\al_s)=\al_s\be_{t,s}^{-1}\al_t^{-m}\be_{t,s} \qquad\text{ (see Figure \ref{Fig9}).}
\]
\endproof

\begin{prop}\label{P5}
Let $x\in\pi_1(F(\Gamma),\basepts)$, and let $s,t\in S$ and $m\in\Z$ such that $|m|\geq 2$. Suppose that $x$ is without squares in $\al_t$ and that $\tau_t^m(x)$ has a square in $\al_s$. Then one of the following two conditions hold:
\begin{description}
\item{(a)} $s=t$;
\item{(b)} $m_{s,t}=2$ and $x$ has a square in $\al_s$.
\end{description}
\end{prop}

\Proof
Let $x=\mu_0\al_{s_1}^{k_1}\mu_1...\al_{s_\ell}^{k_\ell}\mu_\ell$ in reduced form. Since $x$ is without squares in $\al_t$, we have $k_i=\pm 1$ if $s_i=t$. By Proposition \ref{P4}, the reduced form of $\tau_t^m(x)$ may be written
\[\tau_t^m(x)=x_0v_1x_1...v_\ell x_\ell\]
where:
\begin{description}
\item{$\bullet$} $x_i$ is a reduced form of type $(\mu_i,\al_t^m)$ for all $i=1,...,\ell$;
\item{$\bullet$} if $s_i=t$, then $v_i=\al_{s_i}^{p_i}$ where $p_i\equiv\pm 1$ (mod $m$);
\item{$\bullet$} if $m_{s_i,t}=2$, then $v_i=\al_{s_i}^{k_i}$;
\item{$\bullet$} if $m_{s_i,t}=3$, then $v_i= \al_{s_i}^{\varep_i}(v_i'\al_{s_i}^{\varep_i})^{|k_i|-1}$
where $\varep_i\in\{\pm 1\}$ and $v_i'$ is a non-constant reduced form of type $(1,\al_t^m)$.
\end{description}
Note that in the last case, when $m_{s_i,t}=3$, we actually have $\tau_t^m(\al_{s_i}^{k_i})=(v_i'\al_{s_i}^{\varep_i})^{|k_i|}$ or $(\al_{s_i}^{\varep_i}v_i')^{|k_i|}$ depending on whether $s_i<t$ or not and on whether $k_i$ is positive or negative. In each case, the extra $v_i'$ is absorbed into the adjacent form $x_{i-1}$ or $x_i$ respectively. 

It now follows that, if $\tau_t^m(x)$ has a square in $\al_s$, then either $s=t$ or there exists $i\in\{ 1,...,\ell\}$ such that $m_{s_i,t}=2$ and $|k_i|\geq 2$.
\endproof

\section{The small type case}\label{S5}

For the moment, suppose that $\Gamma$ is an arbitrary Coxeter graph. Take an abstract set $\Cal T=\{T_s: s\in S\}$ in bijection with $S$, and consider the group $H(\Gamma)$ defined by the presentation
\[
H(\Gamma)=\< \Cal T\mid T_sT_t=T_tT_s \text{ whenever } m_{s,t}=2 \>\,.
\]
Note that this presentation depends only on the unlabelled Coxeter graph, and that $H(\Gamma)=A(\Gamma_\infty)$ where $\Gamma_\infty$ is the Coxeter graph obtained from $\Gamma$ by relabelling every edge with $\infty$.  The class of groups having such a presentation, known as ``graph groups'' or ``right-angled Artin groups'', have been widely studied. They are examples of so-called ``graph products'' of groups (namely the graph products of infinite cyclic groups) which were studied in the thesis of E.R. Green \cite{Gr} (see also \cite{HM}).
\medskip

Suppose that, for each $s\in S$, we are given an integer $m_s\geq 2$. We then have a homomorphism $f:H(\Gamma)\to A(\Gamma)$ which sends $T_s$ to $\sig_s^{m_s}$ for all $s\in S$. It is our objective in this paper (in proving Theorem \ref{main}) to show that this homomorphism $f$ is always injective. In this section, we shall treat the small type case.
\medskip

We begin with some general notions concerning the presentation given above for the group $H(\Gamma)$.  An \demph{expression} in $\Cal T$ is a sequence
\[
W = (T_{s_\ell}^{p_\ell},...,T_{s_2}^{p_2},T_{s_1}^{p_1})
\]  
where $s_i\in S$, $p_i\in\Z\setminus\{ 0\}$.  The number $\ell$ is called the \demph{length} of the expression $W$, written $\ell=\ell(W)$. We write $[W]$ for the corresponding element $T_{s_\ell}^{p_\ell}...T_{s_2}^{p_2}T_{s_1}^{p_1}\in H(\Gamma)$. Given $w\in H(\Gamma)$, we say that $W$ is \demph{an expression for $w$} if $w=[W]$.
The following terminology is motivated by that used by Brown in \cite{Br}.

Consider the expression $W$ as written above. Suppose that there exists $i\in\{ 1,...,\ell-1\}$ such that $s_i=s_{i+1}$. Put
\[
W'=
\begin{cases}
(..., T_{s_{i+2}}^{p_{i+2}}, T_{s_i}^{p_{i+1}+p_i}, 
T_{s_{i-1}}^{p_{i-1}},...)\qquad 	&\text{if }p_i+p_{i+1}\neq 0,\\
(...,T_{s_{i+2}}^{p_{i+2}}, T_{s_{i-1}}^{p_{i-1}},...)\qquad&\text{if }p_i+p_{i+1}=0.
\end{cases}
\]
We say that $W'$ is obtained from $W$ via an \demph{elementary $M$-operation of type I}. This operation shortens the length of an expression by $1$ or $2$.

Suppose that there exists $i\in\{ 1,...,\ell-1\}$ such that $m_{s_i,s_{i+1}}=2$. Put
\[
W''= (..., T_{s_{i+2}}^{p_{i+2}}, T_{s_i}^{p_i}, T_{s_{i+1}}^{p_{i+1}}, 
T_{s_{i-1}}^{p_{i-1}},...)
\]
We say that $W''$ is obtained from $W$ via an \demph{elementary $M$-operation of type II}. This operation leaves the length of an expression unchanged.

We shall say that $W$ is \demph{$M$-reduced} if the length of $W$ can not be reduced by applying a sequence of elementary $M$-operations. Clearly, every element of $H(\Gamma)$ has an $M$-reduced expression.
If two expressions $W$ and $W'$ are related by a sequence of elementary $M$-operations of type II, then we say that they are \demph{II-equivalent} and write $W\simII W'$.

Let $s\in S$. An $M$-reduced expression is said to \demph{end in $s$} if it is  II-equivalent to an expression $(T_{s_\ell}^{p_\ell},...,T_{s_1}^{p_1})$ in which $s_\ell=s$.

\begin{lemma}\label{Mred}
Suppose that $W=(T_{s_\ell}^{p_\ell},...,T_{s_2}^{p_2},T_{s_1}^{p_1})$ is an expression, and fix an integer $r$ such that $2\leq r\leq \ell$. Then the following are equivalent:
\begin{description}
\item{(i)} $W$ is $M$-reduced;
\item{(ii)} for all $i<j$ such that $s_i=s_j$, there exists $k$ such that $i<k<j$ and $m_{s_i,s_k}\neq 2$;
\item{(iii)} the expressions $U=(T_{s_r}^{-p_r},...,T_{s_\ell}^{-p_\ell})$ and $V=(T_{s_{r-1}}^{p_{r-1}},...,T_{s_1}^{p_1})$ are both $M$-reduced, and there is no $s\in S$ such that both $U$ and $V$ end in $s$.
\end{description}
\end{lemma}

\Proof That both (ii) and (iii) are consequences of (i) is obvious. We note that the truth of condition (ii) is invariant under elementary $M$-operations of type II, while an elementary $M$-operation of type I can never be performed on an expression satisfying (ii). Thus (ii)$\imply$(i). We now show that (iii)$\imply$(ii). Suppose that (iii) holds, but that there exist $i<j$ such that $s_i=s_j$ and $m_{s_i,s_k}=2$ for all $i<k<j$. If $i<j<r$, or $r\leq i<j$, then we contradict that both $V$ and $U$ are $M$-reduced, while if $i<r\leq j$ then it is clear that both $U$ and $V$ end in $s_i$.
\endproof

\medskip

We return now to the case where $\Gamma$ is a connected Coxeter graph of small type. Recall the homomorphism $g:A(\Gamma)\to\Cal M(F(\Gamma))$ of Proposition \ref{P2}. We now consider the group $H(\Gamma)$ to act on $\pi_1(F(\Gamma),\basepts)$ via the homomorphism $g\circ f:H(\Gamma)\to A(\Gamma)\to\Cal M(F(\Gamma))$.

\begin{prop}\label{P6}
Let $W$ be an $M$-reduced expression for $w\in H(\Gamma)$. Let   $x\in\pi_1(F(\Gamma),\basepts)$, and $t\in S$. Suppose that $x$ is without squares in $\al_s$ for all $s\in S$ and that $w(x)$ has a square in $\al_t$. Then $W$ ends in $t$.
\end{prop}

\Proof
We proceed by induction on the length of $W$.
Let $W$ be written $(T_{s_\ell}^{p_\ell},...,T_{s_1}^{p_1})$ and write
$W' = (T_{s_{\ell-1}}^{p_{\ell-1}},...,T_{s_1}^{p_1})$. If $[W'](x)$ has a square in $\al_{s_\ell}$ then, by the induction hypothesis, $W'$ ends in $s_\ell$. But then $W$ is II-equivalent to an expression $(T_{s_\ell}^{p_\ell},T_{t_{\ell-1}}^{q_{\ell-1}},...,T_{t_1}^{q_1})$, where $t_{\ell-1}=s_\ell$. But this expression admits an elementary operation of type I, contradicting the assumption that $W$ is $M$-reduced. Thus, $[W'](x)$ is without squares in $\al_{s_\ell}$.

Note that $w(x)=\tau_{s_\ell}^m([W'](x))$ where $m=p_\ell m_{s_\ell}$ and $|m|\geq 2$. We are given that $w(x)$ has a square in $\al_t$. Therefore, by Proposition \ref{P5}, one either has that $s_\ell =t$ or that $m_{s_\ell,t}=2$ and $[W'](x)$ has a square in $\al_t$. If $s_\ell=t$, then $W$ obviously ends in $t$. We may therefore suppose the latter, in which case, by induction, $W'$ ends in $t$. In other words, $W'$ is II-equivalent to an expression  $(T_{t_{\ell-1}}^{q_{\ell-1}},...,T_{t_1}^{q_1})$ with $t_{\ell -1}=t$. Then, using the fact that $m_{s_\ell,t}=2$, 
\[
W\simII (T_{s_\ell}^{p_\ell},T_t^{q_{\ell -1}}, T_{t_{\ell -2}}^{q_{\ell -2}}, ..., T_{t_1}^{q_1})\simII (T_t^{q_{\ell-1}}, T_{s_\ell}^{p_\ell}, T_{t_{\ell-2}}^{q_{\ell-2}}, ..., T_{t_1}^{q_1})\,.
\]
Hence $W$ ends in $t$.
\endproof

\begin{prop}\label{P7}
If there exists a nontrivial $M$-reduced expression for $w\in H(\Gamma)$,   then $f(w)$ is nontrivial. In particular, the homomorphism $f:H(\Gamma)\to A(\Gamma)$ is injective.
\end{prop}

\Proof
Note that the result is obvious in the case that $\Gamma$ is a trivial graph (i.e: $\Gamma$ consists of just one vertex, and both $H(\Gamma)$ and $A(\Gamma)$ are infinite cyclic). Suppose then that $\Gamma$ is nontrivial, and let $W=(T_{s_\ell}^{p_\ell},...,T_{s_1}^{p_1})$ be a nontrivial $M$-reduced expression for $w\in H(\Gamma)$. Using the fact that $\Gamma$ is connected and nontrivial, choose $t\in S$ such that $m_{s_\ell ,t}=3$. We shall show that $w(\al_t)\neq\al_t$. This clearly implies that $f(w)\neq 1$.

Suppose that $w(\al_t)=\al_t$. Putting $W'=(T_{s_{\ell -1}}^{p_{\ell -1}},...,T_{s_1}^{p_1})$, we have \[ [W'](\al_t) = T_{s_\ell}^{-p_\ell}(\al_t) = \tau_{s_\ell}^m(\al_t) \] where $m=-p_\ell m_{s_\ell}$, and $|m|\geq 2$. By Proposition \ref{P4}, $\tau_{s_\ell}^m(\al_t)$ has a square in $\al_{s_\ell}$. It follows, by Proposition \ref{P6}, that $W'$ ends in $s_\ell$. But this contradicts $W$ being $M$-reduced. Therefore $w(\al_t)\neq\al_t$. 
\endproof

\begin{corollary}
Theorem \ref{main} and the Tits Conjecture hold for irreducible Artin groups of small type.
\end{corollary}

It has been known for some time that the finite process of obtaining an $M$-reduced expression from any given expression (via elementary $M$-operations) provides a solution to the word problem of the graph group $H(\Gamma)$. Namely: 
\begin{equation}\label{wordproblem}
\text{If $W$ is an $M$-reduced expression representing the identity then $W$ is trivial.} 
\end{equation}
We refer to \cite{Gr}, and more recently \cite{HW}, both of which treat the more general case of graph products. We point out that, none of our arguments rely upon this solution to the word problem. On the contrary, statement (\ref{wordproblem}) is actually a \emph{consequence} of Proposition \ref{P7} (for graph groups only). We note also that, by using Lemma \ref{Mred}(iii) and induction on the lengths of expressions, 
one can show the strengthened version of (\ref{wordproblem}), as stated in \cite{HM}, that
any two $M$-reduced expressions for the same element are II-equivalent.

\section{The general case}\label{S6}

\begin{figure}[ht]
\begin{center}
\includegraphics[width=15cm]{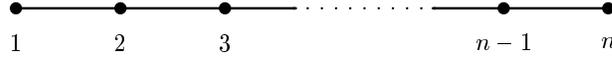}
\end{center}
\caption{The Coxeter graph $A_n$.}\label{Fig10}
\end{figure}

 We introduce the Coxeter graph $A_n$ (for $n>0$) of Figure \ref{Fig10}, and number the vertices of $A_n$ as shown. The Artin group $A(A_n)$ is better known as the braid group on $(n+1)$ strings.

We extend the notation $\pro(-,-;m)$ defined in the introduction to apply to arbitrary expressions $x,y$. Thus  $\pro(x,y;m)= (xy)^{\frac{m}{2}}$ if $m$ is even, and $x(yx)^{\frac{m-1}{2}}$ if $m$ is odd. The following result
can be found in the paper of Brieskorn and Saito \cite{BS}, Lemma 5.8.

\begin{prop}\label{P10}
For $n\in\N\setminus\{ 0\}$, let $\{n,n-1\}=\{k,\ell\}$ where $k$ is odd and $\ell$ is even, and   let $P=\sig_1\sig_3...\sig_k$ and $Q=\sig_2\sig_4...\sig_\ell$ (or $1$ if $n=1$) denote elements of $A(A_n)$. Then 
\[
\pro(P,Q;n+1)=\pro(Q,P;n+1)\,.
\]
\end{prop}

The element $\pro(P,Q;n+1)$ of $A(A_n)$ is the so-called `fundamental element' introduced into the study of braid groups by Garside \cite{Ga}. 

\begin{figure}[ht]
\begin{center}
\includegraphics[width=15cm]{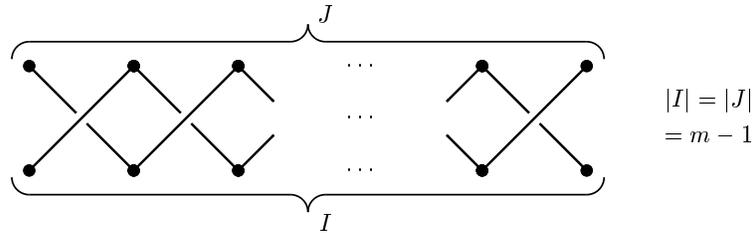}
\end{center}
\caption{The Coxeter graph $\Gamma(m)$.}\label{Fig11}
\end{figure}

Let $m$ be an integer, $m\geq 3$. Let $\Gamma(m)$ denote the Coxeter graph with $2(m-1)$ vertices illustrated in Figure \ref{Fig11}. This  is a bi-partite graph with vertex set $I\sqcup J$, where $|I|=|J|=m-1$ as shown in the figure. As a Coxeter graph, $\Gamma(m)$ is the disjoint union of two copies of $A_{m-1}$. Thus $A(\Gamma(m))\cong A(A_{m-1})\times A(A_{m-1})$.

Let $\al,\be$ be the elements of $A(\Gamma(m))$ defined by
\[
\al=\prod\limits_{i\in I}\sig_i\qquad \text{and}\qquad \be=\prod\limits_{j\in J} \sig_j
\]
Then, by Proposition \ref{P10}, we have
\[\pro(\al,\be ;m)=\pro(\be,\al ;m)\,.\]

Let $k\in\N$. We let $k\Gamma(m)$ denote the disjoint union of $k$ copies of $\Gamma(m)$, which is a bi-partite graph with vertex set $kI\sqcup kJ$, where $kI$ denotes the $k$ disjoint copies of $I$, and $kJ$ the $k$ disjoint copies of $J$. Now let $\al,\be$ be the elements of $A(k\Gamma(m))$ defined by
\[
\al=\prod\limits_{i\in kI}\sig_i\qquad \text{and}\qquad \be=\prod\limits_{j\in kJ} \sig_j
\]
We still have
\[\pro(\al,\be ;m)=\pro(\be,\al ;m)\,.\]

\begin{figure}
\begin{center}
\includegraphics[width=15cm]{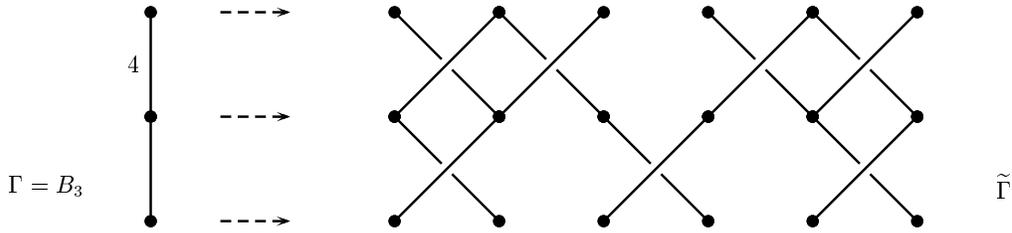}
\end{center}
\caption{A graph $\wtil\Gamma$ for $\Gamma=B_3$.}\label{Fig12}
\end{figure}

We suppose now that $\Gamma$ is a connected Coxeter graph in which $m_{s,t}\neq\infty$ for all $s,t\in S$.
Let $N$ denote the least common multiple of the set $\{ m_{s,t}-1 : s\neq t\in S \}$ of natural numbers. For each $s\in S$, take an abstract set $I(s)$ with $N$ elements. We define a Coxeter graph $\wtil\Gamma$ of small type as follows:
\begin{description}
\item{$\bullet$} the vertex set of $\wtil\Gamma$ is the disjoint union of the sets $I(s)$ for $s\in S$;
\item{$\bullet$} if $m_{s,t}=2$, then there are no edges between the vertices of  $I(s)$ and $I(t)$;
\item{$\bullet$} if $m_{s,t}\geq 3$, then the full subgraph of $\wtil\Gamma$ spanned by $I(s)\sqcup I(t)$ is isomorphic to $(\frac{N}{m_{s,t}-1})\Gamma(m_{s,t})$ by a graph isomorphism taking $I(s)$ to $(\frac{N}{m_{s,t}-1})I$ and $I(t)$ to $(\frac{N}{m_{s,t}-1})J$.
\end{description}
Such a graph always exists, but is not unique in general. An example is given in 
Figure \ref{Fig12}.

 By the preceding remarks, we now have:

\begin{prop}\label{P11}
There exists a homomorphism $\phi : A(\Gamma)\to A(\wtil\Gamma)$ which sends $\sig_s$ to $\prod_{i\in I(s)}\sig_i$ for each $s\in S$.
\end{prop} 

\Remark The homomorphism of Proposition \ref{P11} is an example of a homomorphism obtained by a ``folding'' of $\wtil\Gamma$ onto $\Gamma$ as described in \cite{JC}. It was shown in \cite{JC} that the corresponding homomorphism $\phi^+$ between Artin monoids is injective, and it seems quite possible (though this remains a conjecture) that $\phi$ is also injective.
\medskip

Recall that one is given an integer $m_s\geq 2$ for each $s\in S$ and that $f: H(\Gamma)\to A(\Gamma)$ denotes the homomorphism which sends $T_s$ to $\sig_s^{m_s}$ for each $s\in S$.

\begin{prop}\label{P12}
Let $\Gamma$ be a connected Coxeter graph with $m_{s,t}\neq\infty$ for all $s,t\in S$. Then the homomorphism $f: H(\Gamma)\to A(\Gamma)$ is injective.
\end{prop} 

\Proof
Consider the homomorphism $\wtil f : H(\wtil\Gamma)\to A(\wtil\Gamma)$ which sends $T_i$ to $\sig_i^{m_s}$ for all $s\in S$ and all $i\in I(s)$, and the homomorphism $\psi :  H(\Gamma)\to H(\wtil\Gamma)$ which sends $T_s$ to $\prod_{i\in I(s)}T_i$. Then we clearly have $\wtil f\circ\psi = \phi\circ f$. Since $\wtil f$ is injective by Proposition \ref{P7}, it suffices to show that $\psi$ is injective. 

Let $W = (T_{s_\ell}^{p_\ell},...,T_{s_2}^{p_2},T_{s_1}^{p_1})$ be an $M$-reduced expression for a nontrivial element $w\in H(\Gamma)$. Label the vertices of $\wtil\Gamma$ so that $I(s)=\{s(1),s(2),...,s(N)\}$, for each $s\in S$. Then the expression 
\[
U = (T_{s_\ell(N)}^{p_\ell},..,T_{s_\ell(1)}^{p_\ell},....,T_{s_2(N)}^{p_2},..,
T_{s_2(1)}^{p_2},T_{s_1(N)}^{p_1},..,T_{s_1(1)}^{p_1})
\]
is an expression for $\psi(w)$, and is $M$-reduced by condition (ii) of Lemma \ref{Mred}. Clearly $U$ is nontrivial (since $W$ is), and so, by Proposition \ref{P7}, $\psi(w)\neq 1$.
\endproof

\begin{prop}\label{P13}
Let $\Gamma$ be an arbitrary Coxeter graph.
Then the homomorphism $f: H(\Gamma)\to A(\Gamma)$ is injective.
\end{prop} 

\Proof
First suppose that $\Gamma$ is connected. We define $\what\Gamma$ to be the Coxeter graph associated to the Coxeter matrix $\what M=(\what m_{s,t})_{s,t\in S}$ where
\[
\what m_{s,t}=
\begin{cases}
\hfil 3\quad&\text{if } m_{s,t}=\infty\\
m_{s,t}\quad&\text{otherwise.} 
\end{cases}
\]

We have a homomorphism $\phi :A(\Gamma)\to A(\what\Gamma)$ sending $\sig_s$ to $\sig_s$ for each $s\in S$, and a homomorphism $\what f: H(\Gamma)=H(\what\Gamma)\to A(\what\Gamma)$ sending $T_s$ to $\sig_s^{m_s}$ for each $s\in S$. Since $\what f= \phi\circ f$ and $\what f$ is injective, by Proposition \ref{P12}, it follows that $f$ is injective.

Finally, suppose that $\Gamma$ is a disjoint union of connected subgraphs $\Gamma_i$ for $i=1,..,k$. Then $A(\Gamma)$ and $H(\Gamma)$ are direct products 
\[
A(\Gamma)=\prod\limits_{i=1}^k A(\Gamma_i)\,,\quad\text{and}\quad
H(\Gamma)=\prod\limits_{i=1}^k H(\Gamma_i)\,,
\]
and the homomorphism $f$ is a product $f_1\times f_2\times\cdot\cdot\times f_k$ of injective homomorphisms $f_i:H(\Gamma_i)\to A(\Gamma_i)$. Thus $f$ is injective.
\endproof

This completes the proof of Theorem \ref{main} and the Tits conjecture for an arbitrary Artin group.

\end{document}